\documentclass[final,1p,times]{elsarticle}

\usepackage{amssymb}
 \usepackage{amsthm}
\usepackage{amscd}
\usepackage{amsmath}
\usepackage{amsfonts}
\usepackage{amssymb}
\usepackage{graphicx}
\usepackage{color}
\newtheorem{theorem}{Theorem}
\newtheorem{proposition}[theorem]{Proposition}
\newtheorem{lemma}[theorem]{Lemma}

\usepackage{mathrsfs}
\usepackage{titletoc}

\newcommand{\ra}{\rightarrow}
\newcommand{\p}{\partial}
\newcommand{\f}{\frac}

\renewcommand{\f}{\frac}

\newcommand{\be}{\begin{equation}}
\newcommand{\ee}{\end{equation}}
\newcommand{\bea}{\begin{eqnarray}}
\newcommand{\eea}{\end{eqnarray}}
\newcommand{\bna}{\begin{eqnarray*}}
\newcommand{\ena}{\end{eqnarray*}}

\renewcommand{\le}{\left}
\newcommand{\ri}{\right}

\journal{***}

\begin{document}

\begin{frontmatter}

\title{Nonexistence of extremals for an inequality of Adimurthi-Druet
on a closed Riemann surface}

\author{Yunyan Yang}
 \ead{yunyanyang@ruc.edu.cn}

 \address{Department of Mathematics,
Renmin University of China, Beijing 100872, P. R. China}

\begin{abstract}

Based on a recent work of Mancini-Thizy \cite{Mancini-Thizy}, we obtain the nonexistence of extremals for an inequality
of Adimurthi-Druet \cite{A-D} on a closed Riemann surface $(\Sigma,g)$. Precisely, if $\lambda_1(\Sigma)$ is the first eigenvalue
of the Laplace-Beltrami operator with respect to the zero mean value condition, then there exists a positive real number $\alpha^\ast<\lambda_1(\Sigma)$ such that for all $\alpha\in
(\alpha^\ast,\lambda_1(\Sigma))$, the supremum
 $$\sup_{u\in W^{1,2}(\Sigma,g),\,\int_\Sigma udv_g=0,\,\|\nabla_gu\|_2\leq 1}\int_\Sigma \exp(4\pi u^2(1+\alpha\|u\|_2^2))dv_g$$
 can not be attained by any $u\in W^{1,2}(\Sigma,g)$ with $\int_\Sigma udv_g=0$ and $\|\nabla_gu\|_2\leq 1$, where $W^{1,2}(\Sigma,g)$ denotes
 the usual Sobolev space and $\|\cdot\|_2=(\int_\Sigma|\cdot|^2dv_g)^{1/2}$ denotes the $L^2(\Sigma,g)$-norm.
  This complements our earlier result in \cite{Yang-Tran}.

 \end{abstract}


\end{frontmatter}

\section{Introduction}
Let $\Omega$ be a smooth bounded domain of $\mathbb{R}^n$, $n\geq 2$, and  $W_0^{1,p}(\Omega)$ be the usual
Sobolev space, namely, the completion of $C_0^\infty(\Omega)$ under the norm
$$\|u\|_{W_0^{1,p}(\Omega)}=\le(\int_\Omega|\nabla u|^pdx\ri)^{1/p}$$
for some $p\geq 1$. By the Sobolev embedding theorem, when $p<n$, $W_0^{1,p}(\Omega)$ is continuously embedded in $L^{\f{np}{n-p}}(\Omega)$;
when $p>n$, $W_0^{1,p}(\Omega)$ is continuously embedded in $C^{1-\f{n}{p}}(\Omega)$. However $W_0^{1,n}(\Omega)$ can not be embedded in $L^\infty(\Omega)$.
This limit case was solved respectively by Yudovich \cite{Yudovich}, Pohozaev \cite{Pohozaev}, Peetre \cite{Peetre}, Trudinger \cite{Trudinger} and Moser \cite{Moser}.
Precisely there holds
\be\label{T-M-ineq}\sup_{u\in W_0^{1,n}(\Omega),\,\|u\|_{W_0^{1,n}(\Omega)}\leq 1}\int_\Omega \exp(\alpha_n|u|^{\f{n}{n-1}})dx<\infty,\ee
where $\alpha_n=n\omega_{n-1}^{1/(n-1)}$ and $\omega_{n-1}$ denotes the area of the unit sphere in $\mathbb{R}^n$. Moreover, $\alpha_n$ is the best constant
in the sense that if $\alpha_n$ is replaced by any larger number $\alpha$, the integrals in (\ref{T-M-ineq}) are still finite, but the supremum is
infinite. The existence of extremals for the supremum in (\ref{T-M-ineq}) was first obtained by Carleson-Chang \cite{CC} in the case that $\Omega$ is a unit ball
in $\mathbb{R}^n$, then extended by Struwe \cite{Struwe} to a domain close to a unit ball in the sense of measure, by Flucher \cite{Flucher} to
a general domain in $\mathbb{R}^2$, and by Lin \cite{Lin} to an arbitrary domain in $\mathbb{R}^n$.

In literature, (\ref{T-M-ineq}) is known as the Trudinger-Moser inequality. It was improved by Adimurthi-Druet \cite{A-D} as follows:
Let $\Omega$ be a smooth bounded domain of $\mathbb{R}^2$,
and $\lambda_1(\Omega)$ be the first eigenvalue of the Laplace-Beltrami operator with respect to the Dirichlet boundary condition. Then for any
$\alpha<\lambda_1(\Omega)$, there  holds
\be\label{adi-druet}\sup_{u\in W_0^{1,2}(\Omega),\,\|u\|_{W_0^{1,2}(\Omega)}\leq 1}\int_\Omega \exp(4\pi u^2(1+\alpha\|u\|_2^2))dx<\infty;\ee
moreover, the above supremum is infinite if $\alpha\geq\lambda_1(\Omega)$, where $\|\cdot\|_p=(\int_\Omega|\cdot|^pdx)^{1/p}$.
This result was generalized by the author \cite{JFA-06} to the higher dimensional case. In particular,
$\Omega$ is assumed to be a smooth bounded domain in $\mathbb{R}^n$, $n\geq 3$. Define
$\lambda_n(\Omega)=\inf_{u\in W_0^{1,n}(\Omega),\,u\not\equiv 0}\|\nabla u\|_n^n/\|u\|_n^n$. Then for any $\alpha<\lambda_n(\Omega)$,
there holds
\be\label{high-dimension}\sup_{u\in W_0^{1,n}(\Omega),\,\|u\|_{W_0^{1,n}(\Omega)}\leq 1}\int_\Omega
\exp(\alpha_n|u|^{\f{n}{n-1}}(1+\alpha\|u\|_n^n)^{\f{1}{n-1}})dx<\infty,\ee
where $\alpha_n$ is defined as in (\ref{T-M-ineq}). Moreover, the supremum in (\ref{high-dimension}) can be attained by some
$u\in W_0^{1,n}(\Omega)$ with $\|u\|_{W_0^{1,n}(\Omega)}=1$ for all $\alpha<\lambda_n(\Omega)$.
But in the case $n=2$, it was shown by Lu-Yang \cite{Lu-Yang-09} that extremals of the supremum in (\ref{adi-druet}) exist only for sufficiently
small $\alpha>0$.
However, among other results, an inequality stronger than
(\ref{adi-druet}) was obtained by Tintarev \cite{Tintarev}. Namely, if $\alpha<\lambda_1(\mathbb{B}_R)$ with
$|\mathbb{B}_R|=|\Omega|$, where $|\cdot|$ denotes the Lebesgue measure of a set, then
\be\label{Tin}\sup_{u\in W_0^{1,2}(\Omega),\, \|\nabla u\|_2^2-\alpha\|u\|_2^2\leq 1}\int_\Omega \exp(4\pi u^2)dx<\infty.\ee
Later, it was shown by the author \cite{JDE-15} that (\ref{Tin}) holds for all $\alpha<\lambda_1(\Omega)$, that
extremals of the supremum in (\ref{Tin}) exist for all $\alpha<\lambda_1(\Omega)$,
and that similar results are still valid when higher order eigenvalues of
the Laplace-Beltrami operator are taken into account. Recently a higher dimensional version of (\ref{Tin}) was established
by Nguyen \cite{Nguyen}, say for any $\alpha<\lambda_n(\Omega)$,
\be\label{Nguy}\sup_{u\in W_0^{1,n}(\Omega),\, \|\nabla u\|_n^n-\alpha\|u\|_n^n\leq 1}\int_\Omega \exp(\alpha_n |u|^{\f{n}{n-1}})dx<\infty;\ee
also the corresponding extremals exist. One can check that (\ref{Nguy}) is stronger than (\ref{high-dimension}).
Nevertheless one may ask whether or not extremals of (\ref{adi-druet}) exist for $\alpha$ sufficiently
close to $\lambda_1(\Omega)$.
Based on works of Malchiodi-Martinazzi \cite{M-Martinazzi}, Mancini-Martinazzi \cite{Mancini-Martinazzi}
and Druet-Thizy \cite{Druet-Thizy}, it was Mancini-Thizy \cite{Mancini-Thizy} who gave a negative answer. Namely, when $\alpha$ is
sufficiently close to $\lambda_1(\Omega)$, the supremum in (\ref{adi-druet}) has no extremal.

Trudinger-Moser inequalities were introduced on Riemannian manifolds by Aubin \cite{Aubin10}, Cherrier \cite{Cherrier1} and
Fontana \cite{Fontana}. Let $(M,g)$ be an $n$-dimensional closed Riemannian manifold, $\alpha_n$ be as in (\ref{T-M-ineq})
and $W^{1,n}(M,g)$ be the usual
Sobolev space. Then there holds
\be\label{Fontana-ineq}
\sup_{u\in W^{1,n}(M),\,\int_Mudv_g=0,\,\int_M|\nabla_gu|^ndv_g\leq 1}\int_M\exp(\alpha_n|u|^{\f{n}{n-1}})dv_g<\infty,
\ee
where $\nabla_g$ denotes the gradient operator and $dv_g$ stands for the volume element.
Based on works of Ding-Jost-Li-Wang \cite{DJLW} and Adimurthi-Struwe \cite{A-S}, Li \cite{LiJPDE,LiSci} was able to prove the existence
of extremals for the supremum in (\ref{Fontana-ineq}). In \cite{Yang-Tran}, the author improved (\ref{Fontana-ineq}) in the case $n=2$
as follows:
Assuming that $(\Sigma,g)$ is a closed Riemann surface, we get an analog of (\ref{adi-druet}), namely
 \be\label{A-D-ineq}\sup_{u\in W^{1,2}(\Sigma),\,\int_\Sigma udv_g=0,\,\int_\Sigma |\nabla_gu|^2dv_g\leq 1}
 \int_\Sigma \exp(4\pi u^2(1+\alpha\|u\|_2^2))dv_g<\infty\ee
  for any $\alpha<\lambda_1(\Sigma)=\inf_{u\in W^{1,2}(\Sigma,g), \int_\Sigma udv_g=0, \|u\|_2=1}
\|\nabla_gu\|_2^2$; moreover the above supremum is infinite when $\alpha\geq \lambda_1(\Sigma)$.
Furthermore there exists some $\alpha_\ast>0$ such that extremals exist when $\alpha<\alpha_\ast$. Concerning (\ref{Tin}), we also have its analog
\cite{JDE-15}, say for any $\alpha<\lambda_1(\Sigma)$,
\be\label{Tin-mfd}\sup_{u\in W^{1,2}(\Sigma),\,\int_\Sigma udv_g=0,\,\int_\Sigma |\nabla_gu|^2dv_g-\alpha\int_\Sigma u^2dv_g\leq 1}
 \int_\Sigma \exp(4\pi u^2)dv_g<\infty;\ee
 in addition, extremals exist for any $\alpha<\lambda_1(\Sigma)$. Also (\ref{Tin-mfd}) is stronger than (\ref{A-D-ineq}). In view of \cite{Mancini-Thizy}, we suspect that extremals of the supremum in
 (\ref{A-D-ineq}) do not exist for $\alpha$ sufficiently close to $\lambda_1(\Sigma)$. Our aim is to confirm this suspicion.
 Define a function space
\be\label{H}\mathcal{H}=\le\{u\in W^{1,2}(\Sigma,g): \int_\Sigma udv_g=0\ri\}.\ee
 Then $\lambda_1(\Sigma)$ can  be equivalently written as
 \be\label{eigen}\lambda_1(\Sigma)=\inf_{u\in\mathcal{H},\,u\not\equiv 0}{\|\nabla_gu\|_2^2}/{\|u\|_2^2},\ee
 where $\|\cdot\|_p=(\int_\Sigma|\cdot|^pdv_g)^{1/p}$. We denote for simplicity
\be\label{sup}\Lambda_\alpha(\Sigma)=\sup_{u\in\mathcal{H},\,\|\nabla_gu\|_2\leq 1}\int_\Sigma e^{4\pi u^2(1+\alpha\|u\|_2^2)}dv_g.\ee

 Our main result reads
\begin{theorem}\label{theorem-1}
Let $(\Sigma,g)$ be a closed Riemann surface, $\lambda_1(\Sigma)$ and $\Lambda_\alpha(\Sigma)$ be defined as in (\ref{eigen}) and
(\ref{sup}) respectively. There exists some $\alpha^\ast<\lambda_1(\Sigma)$, such that if $\alpha^\ast<\alpha<\lambda_1(\Sigma)$,
then the supremum $\Lambda_\alpha$ has no extremal.
\end{theorem}

The proof of Theorem \ref{theorem-1} is essentially based on the method of energy estimate, which was used by Mancini-Martinazzi
\cite{Mancini-Martinazzi} and Mancini-Thizy \cite{Mancini-Thizy}.
Let us describe its outline. Suppose Theorem \ref{theorem-1} dose not hold. There would exist a sequence of numbers $\alpha_k$ increasingly
tending to $\lambda_1(\Sigma)$ such that $\Lambda_{\alpha_k}=\int_\Sigma e^{4\pi u_k^2(1+\alpha_k\|u_k\|_2^2)}dv_g$ for some function $u_k\in\mathcal{H}$
with $\|\nabla_gu_k\|_2=1$. Clearly $u_k$ satisfies the corresponding Euler-Lagrange equation (see (\ref{E-L}) below).
Moreover $u_k\ra 0$ strongly in $L^p(\Sigma,g)$
for any $p>0$.
By performing local blow-up analysis
and global analysis on $u_k$, we obtain
$$\int_\Sigma|\nabla_gu_k|^2dv_g=1-\f{3}{2}\alpha_k^2\|u_k\|_2^4+o(\|u_k\|_2^4),$$
which contradicts $\|\nabla_gu_k\|_2=1$ for sufficiently large $k$, since $\alpha_k\ra\lambda_1(\Sigma)$. It should be remarked that
since $u_k$ changes sign in our case, we must re-establish a gradient estimate on $u_k$ instead of
 Druet's original one \cite{Druet}; moreover, when performing local blow-up analysis, we must take into account
the fact that $u_k$ may change sign near blow-up points. Here we choose a sequence of isothermal coordinates instead of a fixed isothermal coordinate.
Such coordinates greatly simplify the energy calculation.
    When estimating the energy on domains away from blow-up points,
we use the H\"older inequality and the global convergence of $u_k/\|u_k\|_p$ for any $p\geq 2$. This is different from (\cite{Mancini-Thizy}, Step 3.4).

Before ending this introduction, we mention related works such as de Souza-do \'O \cite{doo2,doo3}, Ishiwata \cite{Ishiwata}, Martinazzi \cite{Mart}, Martinazzi-Struwe \cite{Mar-Stru},
Lamm-Robert-Struwe \cite{L-R-S}, Adimurthi-Yang \cite{Adi-Yang-Calc}, del Pino-Musso-Ruf \cite{d-M-Ruf1,d-M-Ruf2}, Yang \cite{Yang-Archiv}
and Figueroa-Musso \cite{Figueroa-Musso}. The remaining part of this paper will be organized as follows: In Section \ref{Sec2},
we prove Theorem \ref{theorem-1} by using an energy estimate (Proposition \ref{proposition}); In Section \ref{Sec3}, we prove Proposition
\ref{proposition} by using blow-up analysis. Hereafter we do  not distinguish sequence and subsequence; moreover, we often denote
various constants by the same $C$.

\section{Proof of Theorem \ref{theorem-1}}\label{Sec2}
Let $\mathcal{H}$, $\lambda_1(\Sigma)$ and $\Lambda_\alpha(\Sigma)$ be defined as in (\ref{H}), (\ref{eigen}) and (\ref{sup}) respectively.
In this section, we shall prove Theorem \ref{theorem-1} by contradiction. Suppose the contrary. There would be
a sequence of numbers $(\alpha_k)$ increasingly converging to $\lambda_1(\Sigma)$, such that $\Lambda_{\alpha_k}(\Sigma)$ can be
attained by some function $u_k\in\mathcal{H}$ with $\|\nabla_gu_k\|_2=1$, namely
\be\label{subcritical}\Lambda_{\alpha_k}(\Sigma)=\int_\Sigma e^{4\pi u_k^2(1+\alpha_k\|u_k\|_2^2)}dv_g.\ee
Obviously $(\Lambda_{\alpha_k}(\Sigma))$ is an increasing sequence with respect to $k$. Since the supremum in (\ref{A-D-ineq}) is infinite when
$\alpha\geq \lambda_1(\Sigma)$,  there holds
\be\label{Lambda}\lim_{k\ra\infty}\Lambda_{\alpha_k}(\Sigma)=\Lambda_{\lambda_1(\Sigma)}(\Sigma)=+\infty.\ee
 By a straightforward calculation, $u_k$ satisfies the following Euler-Lagrange equation:
\be\label{E-L}\le\{
\begin{array}{lll}
\Delta_gu_k=\f{\beta_k}{\lambda_k}u_k e^{\sigma_ku_k^2}+\gamma_k u_k-\mu_k \,\,\,{\rm in}\,\,\,\Sigma\\[1.2ex]
u_k\in\mathcal{H},\,\,\, \|\nabla_gu_k\|_2=1\\[1.2ex]
\beta_k=\f{1+\alpha_k\|u_k\|_2^2}{1+2\alpha_k\|u_k\|_2^2},\,\,\, \lambda_k=\int_\Sigma u_k^2e^{\sigma_ku_k^2}dv_g\\[1.2ex]
\sigma_k=4\pi(1+\alpha_k\|u_k\|_2^2),\,\,\, \gamma_k=\f{\alpha_k}{1+2\alpha_k\|u_k\|_2^2}\\[1.2ex]
\mu_k=\f{\beta_k}{{\rm Vol}_g(\Sigma)}\f{1}{\lambda_k}\int_\Sigma u_ke^{\sigma_k u_k^2}dv_g,
\end{array}
\ri.\ee
where $\nabla_g$ and $\Delta_g$ represent the gradient and the Laplace-Beltrami operator respectively. Denote $c_k=\max_{\Sigma}|u_k|$. Since $-u_k$ is also a maximizer of
 $\Lambda_{\alpha_k}(\Sigma)$, in view of (\ref{Lambda}), we can assume up to a subsequence,
 \be\label{c-infty}c_k=u_k(x_k)=\max_\Sigma|u_k|\ra +\infty,\quad x_k\ra x_0\in\Sigma\ee
 as $k\ra \infty$. Note that
 $$\int_\Sigma u_k^2e^{\sigma_ku_k^2}dv_g\geq \int_\Sigma (e^{\sigma_k u_k^2}-1)dv_g,$$
 which together with (\ref{subcritical}) and (\ref{Lambda}) leads to
 \be\label{la-infty}\lambda_k\ra+\infty\quad{\rm as}\quad k\ra\infty.\ee
 To proceed, we observe an energy concentration phenomenon of $u_k$, namely
\begin{lemma}\label{concentration}
 $u_k$ converges to $0$ weakly in $W^{1,2}(\Sigma,g)$, strongly in $L^p(\Sigma)$ for any $p\geq 1$,
and almost everywhere in $\Sigma$. Moreover, $|\nabla_gu_k|^2dv_g\rightharpoonup \delta_{x_0}$
in the sense of measure. As a consequence, there holds
$\sigma_k\ra 4\pi$, $\beta_k\ra 1$ and $\gamma_k\ra\lambda_1(\Sigma)$ as $k\ra\infty$.
\end{lemma}
Since the proof of Lemma \ref{concentration} is an obvious analog of that of (\cite{Yang-Tran}, Lemma 4.3), we omit it,
but refer the reader to \cite{Yang-Tran} for details.
In view of Lemma \ref{concentration}, we have by applying elliptic estimates to (\ref{E-L}) that
\be\label{uk-0}u_k\ra 0\quad{\rm in}\quad C^1_{\rm loc}(\Sigma\setminus\{x_0\}).\ee

We now state the following energy estimate:
\begin{proposition}\label{proposition} Let $(u_k)$ satisfy (\ref{subcritical}) and particularly satisfy (\ref{E-L}).
 Then we have up to a subsequence,
\be\label{Prop-3}\int_\Sigma|\nabla_gu_k|^2dv_g=\f{1}{1+\alpha_k\|u_k\|_2^2}+\f{\alpha_k\|u_k\|_2^4}
{1+2\alpha_k\|u_k\|_2^2}+o(\|u_k\|_2^4),\ee
where $o(\|u_k\|_2^4)/\|u_k\|_2^4\ra 0$ as $k\ra\infty$.
\end{proposition}
The proof of Proposition \ref{proposition} will be postponed to the subsequent section. Assuming this,
we conclude Theorem \ref{theorem-1} as follows:\\

\noindent{\it Completion of the proof of Theorem \ref{theorem-1}.} It follows from Lemma \ref{concentration} that $\|u_k\|_2\ra 0$ as $k\ra\infty$.
Keeping in mind $\|\nabla_gu_k\|_2=1$, we have by (\ref{Prop-3})
 that
\bna
1&=&\int_\Sigma|\nabla_gu_k|^2dv_g\\
&=&1-\alpha_k\|u_k\|_2^2+\f{\alpha_k^2}{2}\|u_k\|_2^4
+\alpha_k\|u_k\|_2^2(1-2\alpha_k\|u_k\|_2^2)+o(\|u_k\|_2^4)\\
&=&1-\f{3}{2}\alpha_k^2\|u_k\|_2^4+o(\|u_k\|_2^4),
\ena
which is a contradiction since $\alpha_k\ra\lambda_1(\Sigma)>0$ as $k\ra\infty$. $\hfill\Box$

\section{Energy estimate}\label{Sec3}

In this section, we prove Proposition \ref{proposition} by analyzing the global and local asymptotic behavior of $u_k$.

\subsection{Isothermal coordinates}
We begin with
the choice of a sequence of
isothermal coordinate systems near blow-up points.
 It is well known (see for example Bers \cite{Bers}, Lecture 3) that there exists an isothermal coordinate system near any point of a Riemann surface.
 Let $(\Sigma,g)$ be the Riemann surface given as in Theorem \ref{theorem-1}.
 In particular, for  $x_0$ given by (\ref{c-infty}), there exists a real number $\delta>0$ and an isothermal coordinate system
 \be\label{isother} (B_{x_0}(\delta),\phi;\{x^1,x^2\}),\ee
 where $B_{x_0}(\delta)\subset\Sigma$ denotes a geodesic ball centered at $x_0$ with radius $\delta$,
  $\phi(x_0)=(0,0)$. In this coordinate system, the metric $g$ can be represented by
\be\label{metr-0}\widetilde{g}(x)=(\phi^{-1})^\ast g(x)=\exp(f(x))(d{x^1}^2+d{x^2}^2),\ee
where $x=(x^1,x^2)$, $f$ is smooth on $\phi(B_{x_0}(\delta))\subset\mathbb{R}^2$, $f(0)=0$ and
\be\label{asy-g}\exp{(f(x))}=1+O(|x|)=1+O\le({\rm dist}_g(\phi^{-1}(x),x_0)\ri).\ee
Here and in the sequel, ${\rm dist}_g(\cdot,\cdot)$ stands for the geodesic distance between two points of $\Sigma$.
Moreover, in the above coordinate system, the gradient operator and
the Laplace-Beltrami operator read as
\be\label{grad-lap}\nabla_{\widetilde{g}}=\exp{(-f)}\nabla_{\mathbb{R}^2}, \quad\Delta_{\widetilde{g}}=-\exp{(-f)}\Delta_{\mathbb{R}^2},\ee
where $\nabla_{\mathbb{R}^2}$ and $\Delta_{\mathbb{R}^2}$ denote the standard gradient operator and Laplacian operator on
$\mathbb{R}^2$ respectively. Recall  $c_k=u_k(x_k)=\max_\Sigma|u_k|\ra+\infty$ and $x_k\ra x_0$ as $k\ra\infty$
(see (\ref{c-infty}) above). Though in many cases of the current topic (see for examples \cite{LiJPDE,Yang-Tran,Yang-Calc,Figueroa-Musso})
the isothermal coordinates (\ref{isother}), especially its properties (\ref{asy-g}) and (\ref{grad-lap}),
had been well used, sometimes (see Section \ref{Sec3.4} below) a sequence of isothermal coordinates will be essentially needed,
namely
\begin{lemma}\label{Lemma-k} There exists some integer $k_0$ such that for any $k\geq k_0$,
one can find an isothermal coordinate system
$(B_{x_k}(\delta/2),\phi_k;\{y^1,y^2\})$ satisfying $\phi_k(x_k)=0$ and
\be\label{g-k}g_k(y)=(\phi_k^{-1})^\ast g(y)=\exp{(f_k(y))}(d{y^1}^2+d{y^2}^2),\ee
where $f_k:\phi_k(B_{x_k}(\delta/2))\ra\mathbb{R}$ is a smooth function with $f_k(0)=0$, $|\nabla^2_{\mathbb{R}^2}f_k|\leq C$,
\be\label{asy-g-k}\exp{(f_k(y))}=1+O(|y|)=1+O\le({\rm dist}_g(\phi_k^{-1}(y),x_k)\ri),\ee
$C^{-1}|y|\leq{\rm dist}_g(\phi_k^{-1}(y),x_k\leq C|y|$ and $|O(|y|)|\leq C|y|$ for some constant $C$ independent of $k$. Moreover, there holds
\be\label{lap-nab}\nabla_{{g_k}}=e^{-f_k}\nabla_{\mathbb{R}^2},\quad\Delta_{{g_k}}=-e^{-f_k}\Delta_{\mathbb{R}^2}.\ee
\end{lemma}

\noindent{\it Proof.} Near the point $x_0$, we first take an isothermal coordinate system $(B_{x_0}(\delta),\phi;\{x^1,x^2\})$
given as in (\ref{isother}). Since $x_k\ra x_0$ as $k\ra\infty$, there must be an integer $k_0$ such that if $k\geq k_0$, then
$B_{x_k}(\delta/2)\subset B_{x_0}(\delta)$.
Thus we have the coordinates $(B_{x_k}(\delta/2),\phi;\{x^1,x^2\})$, in which $\widetilde{x}_k=\phi(x_k)$ and the metric $g$ can be represented
by (\ref{metr-0}). Now we take another coordinate system around $x_k$, namely
$(B_{x_k}(\delta/2),\phi_k;\{y^1,y^2\})$ satisfying
$$y=\phi_k\circ \phi^{-1}(x)=e^{\f{f(\widetilde{x}_k)}{2}}(x-\widetilde{x}_k).$$
Set $f_k(y)=f(x)-f(\widetilde{x}_k)$ and $g_k=(\phi_k^{-1})^\ast g$.
Clearly $\phi_k(x_k)=0$; $f_k:\phi_k(B_{x_k}(\delta/2))\ra\mathbb{R}$ is smooth,  $f_k(0)=0$, and $|\nabla^2_{\mathbb{R}^2}f_k|$ is uniformly bounded
on $\phi_k(B_{x_k}(\delta/2))$, in particular, (\ref{asy-g-k}) holds.
Moreover, since $dy=\exp(f(\widetilde{x}_k)/2)dx$, we have
\bna
g_k(y)&=&\widetilde{g}(x)\\
&=&\exp(f(x))(d{x^1}^2+d{x^2}^2)\\
&=&\exp(f(x))\le(\exp(-f(\widetilde{x}_k))(d{y^1}^2+d{y^2}^2)\ri)\\
&=&\exp(f_k(y))(d{y^1}^2+d{y^2}^2).
\ena
This gives (\ref{g-k}). Finally (\ref{lap-nab}) follows from (\ref{g-k}) immediately.
$\hfill\Box$

\subsection{Global analysis}

From now on, in the above isothermal coordinates $(B_{x_k}(\delta/2), \phi_k; \{y^1,y^2\})$, we write
$\widetilde{x}=\phi_k(x)$ for all $x\in B_{x_k}(\delta/2)$ and $\widetilde{u}=u\circ \phi_k^{-1}$ for all functions $u:B_{x_k}(\delta/2)\ra\mathbb{R}$. Let
\be\label{scale}r_k=\f{\sqrt{\lambda_k}}{\sqrt{\beta_k}\,c_k}\exp{(-\f{\sigma_k}{2}c_k^2)},\ee
where $\lambda_k$, $\beta_k$ and $\sigma_k$ are defined as in (\ref{E-L}).
Using the same argument as in the proof of (\cite{Yang-Tran}, Lemma 4.4), we have that for any fixed real number $a<4\pi$,
\be\label{conv}r_k^2\exp{(a c_k^2)}\ra 0.\ee
It follows from (\ref{E-L}), elliptic estimates and a result of Chen-Li \cite{CL} that
\be\label{bubble}c_k(\widetilde{u}_k(r_ky)-c_k)\ra \varphi(y)=-\f{1}{4\pi}\log(1+\pi|y|^2)\quad
{\rm in}\quad C^2_{\rm loc}(\mathbb{R}^2).\ee
For any $\theta\in(0,1)$, we recall the truncation $u_{k,\theta}=\min\{u_k,\theta c_k\}$ defined
by \cite{LiJPDE,Yang-Tran}. An obvious analog of (\cite{Yang-Tran}, Lemma 4.5) gives that
\be\label{theta}\|\nabla_gu_{k,\theta}\|_2^2=\theta+o_k(1), \,\,\,\|\nabla_g(u_k-u_{k,\theta})\|_2^2=1-\theta+o_k(1).\ee
To understand the asymptotic behavior of $u_k$ away from $x_0$, we need the following:
\begin{lemma}\label{Lemma-5}
$(i)$ $c_k^2/\lambda_k=o_k(1)$; $(ii)$ $\int_\Sigma \f{\beta_k}{\lambda_k}c_k|u_k|\exp{(\sigma_k u_k^2)}dv_g=1+o_k(1)$,
$\int_\Sigma \f{\beta_k}{\lambda_k}c_ku_k\exp{(\sigma_k u_k^2)}dv_g=1+o_k(1)$.
\end{lemma}

\noindent{\it Proof.} Let $\theta\in(0,1)$ be fixed. In view of (\ref{theta}) and the fact that
$\int_\Sigma u_{k,\theta}dv_g=o_k(1)$,
 Fontana's inequality (\ref{Fontana-ineq}) implies that $\exp{(\sigma_ku_{k,\theta}^2)}$ is bounded in $L^p(\Sigma,g)$ for
some $p>1$. This together with Lemma \ref{concentration} leads to
\be\label{L1}\int_{u_k\leq \theta c_k} \exp{(\sigma_ku_{k}^2)}dv_g={\rm vol}_g(\Sigma)+o_k(1).\ee
It then follows that
\bna\int_\Sigma \exp{(\sigma_ku_k^2)}dv_g&=&\int_{u_k>\theta c_k}\exp{(\sigma_ku_k^2)}dv_g+\int_{u_k\leq\theta c_k}\exp{(\sigma_ku_k^2)}dv_g\\
&=&\int_{u_k>\theta c_k}\exp{(\sigma_ku_k^2)}dv_g+{\rm vol}_g(\Sigma)+o_k(1).
\ena
Since $\exp{(\sigma_ku_{k,\theta}^2)}$ is bounded in $L^p(\Sigma,g)$ and $u_k$ converges to $0$ in $L^q(\Sigma,g)$ for any fixed $q>1$, we conclude
by using the H\"older inequality that
\be\label{(2)}\int_{u_k\leq \theta c_k}u_k^2\exp{(\sigma_ku_{k}^2)}dv_g=o_k(1).\ee
Combining (\ref{L1}), (\ref{(2)}) and the definition of $\lambda_k$ (see (\ref{E-L})), we have
\bna\lambda_k&=&\int_{u_k>\theta c_k}u_k^2\exp{(\sigma_ku_k^2)}dv_g+o_k(1)\\
&\geq& \theta^2 c_k^2\le(\int_\Sigma \exp{(\sigma_ku_k^2)}dv_g-{\rm vol}_g(\Sigma)+o_k(1)\ri)+o_k(1).\ena
This together with (\ref{Lambda}) concludes $(i)$.

We now prove $(ii)$. Given any $\theta\in(0,1)$. It follows from $(i)$ and (\ref{L1}) that
 $$\int_{u_k\leq\theta c_k}\f{\beta_k}{\lambda_k}c_k|u_k|\exp{(\sigma_ku_k^2)}dv_g\leq \f{\beta_k c_k^2}{\lambda_k}
 \int_{u_k\leq\theta c_k}\exp{(\sigma_ku_k^2)}dv_g=o_k(1).$$
As a consequence,
\bna
\int_\Sigma\f{\beta_k}{\lambda_k}c_k|u_k|\exp{(\sigma_ku_k^2)}dv_g&=&\int_{u_k>\theta c_k}\f{\beta_k}{\lambda_k}c_ku_k
\exp{(\sigma_ku_k^2)}dv_g+o_k(1)\\
&\leq&\f{1}{\theta}\int_{u_k>\theta c_k}\f{\beta_k}{\lambda_k}u_k^2\exp{(\sigma_ku_k^2)}dv_g+o_k(1)\\
&=&\f{1}{\theta}+o_k(1),
\ena
which leads to
\be\label{leq-1}\limsup_{k\ra\infty}\int_\Sigma\f{\beta_k}{\lambda_k}c_k|u_k|\exp{(\sigma_ku_k^2)}dv_g\leq \f{1}{\theta}.\ee
Nevertheless, by the definition of $\lambda_k$, it is obvious to see
\be\label{leq-2}\liminf_{k\ra\infty}\int_\Sigma\f{\beta_k}{\lambda_k}c_k|u_k|\exp{(\sigma_ku_k^2)}dv_g\geq 1.\ee
Since $\theta\in(0,1)$ is arbitrary, we have by combining (\ref{leq-1}) and (\ref{leq-2}) that
$$\lim_{k\ra\infty}\int_\Sigma\f{\beta_k}{\lambda_k}c_k|u_k|\exp{(\sigma_ku_k^2)}dv_g=1.$$
The other equality in $(ii)$ can be derived in the same way. $\hfill\Box$

\begin{lemma}\label{Lambda-k} There holds
${(\log\lambda_k)}/{c_k^2}=o_k(1)$.
\end{lemma}

\noindent {\it Proof.} Given any $\nu$, $0<\nu<4\pi$. Since $\|\nabla_gu_k\|_2=1$, we have by Fontana's inequality (\ref{Fontana-ineq}) that
$$\int_\Sigma \exp((\sigma_k-\nu)u_k^2)dv_g\leq C$$
for  some constant $C>0$ depending on $\nu$. Hence
$$\lambda_k\leq c_k^2\exp(\nu c_k^2)\int_\Sigma \exp((\sigma_k-\nu)u_k^2)dv_g\leq C c_k^2\exp(\nu c_k^2),$$
which leads to
$$\f{\log\lambda_k}{c_k^2}\leq \nu+\f{2\log c_k}{c_k^2}+\f{C}{c_k^2}.$$
This together with (\ref{c-infty}) and (\ref{la-infty}) implies that
$$0\leq\limsup_{k\ra\infty}\f{\log\lambda_k}{c_k^2}\leq \nu.$$ Since $0<\nu<4\pi$ is arbitrary,
we get the desired result. $\hfill\Box$\\

Moreover, we have the following:
\begin{lemma}\label{Lemma-6} For any $p\geq 2$, there holds
$c_k\|u_k\|_p\ra\infty$ as $k\ra\infty$.
\end{lemma}

\noindent{\it Proof.} By the H\"older inequality, it suffices to prove $c_k\|u_k\|_2\ra\infty$. Suppose not.
There would exist some constant
$C$ such that $c_k\|u_k\|_2\leq C$. Hence $4\pi\alpha_ku_k^2\|u_k\|_2^2\leq C$ on $\Sigma$. This together with
Fontana's inequality (\ref{Fontana-ineq}) leads to
$$\int_\Sigma \exp{(4\pi u_k^2(1+\alpha_k\|u_k\|_2^2))}dx\leq C\int_\Sigma \exp{(4\pi u_k^2)}dv_g\leq C,$$
which contradicts (\ref{Lambda}). $\hfill\Box$\\

Furthermore, we describe the global convergence of $u_k$ as follows:

\begin{lemma}\label{out} For any $p\geq 2$ and any $r>1$, there holds up to a subsequence
${u_k}/{\|u_k\|_p}\ra\psi_p$ in $L^r(\Sigma,g)$, where $\psi_p$ is a smooth solution of the equation
\be\label{ei-func}\le\{\begin{array}{lll}\Delta_g\psi=\lambda_1(\Sigma)\psi\\[1.2ex]
\|\psi\|_p=1.
\end{array}\ri.\ee
\end{lemma}

\noindent{\it Proof.} By (\ref{E-L}), we  have
\be\label{cu}\Delta_g(c_ku_k)=h_k=\f{\beta_k}{\lambda_k}c_ku_k e^{\sigma_ku_k^2}+\gamma_kc_ku_k-\mu_kc_k.\ee
In view of Lemma \ref{Lemma-5}, both $\mu_kc_k$ and $\|\f{\beta_k}{\lambda_k}c_ku_k \exp{(\sigma_ku_k^2)}\|_{1}$ are
bounded. Since $\gamma_k$ is also bounded, we conclude from Lemma \ref{Lemma-6}
that $h_k/(c_k\|u_k\|_p)$ is bounded
in $L^1(\Sigma,g)$ for any $p\geq 2$. In view of (\ref{cu}), one has
\be\label{pp}\Delta_g\f{u_k}{\|u_k\|_p}=\Delta_g\f{c_ku_k}{c_k\|u_k\|_p}=\f{h_k}{c_k\|u_k\|_p}.\ee
For any fixed $r>1$, we take $q<2$ such that $r<2q/(2-q)$. In view of (\ref{pp}),
using the Green representation formula as in the proof of (\cite{Yang-Zhu-Science}, Lemma 2.10),
we have that $u_k/\|u_k\|_p$ is bounded in $W^{1,q}(\Sigma,g)$. Hence there exists some function $\psi_p$
such that $u_k/\|u_k\|_p$ converges to $\psi_p$ weakly in $W^{1,q}(\Sigma,g)$, and strongly in $L^{r}(\Sigma,g)$.
Clearly $\psi_p$ is a distributional solution of (\ref{ei-func}). Applying the regularity theory to (\ref{ei-func}),
we have that $\psi_p$ is smooth. This ends the proof of the lemma. $\hfill\Box$\\

A comparison between $\|u_k\|_p$ and $\|u_k\|_2$ reads

\begin{lemma}\label{2-p}
For any $p\geq 2$, we have $\|u_k\|_p^2=\|u_k\|_2^2/(\|\psi_p\|_2^2+o_k(1))$.
\end{lemma}
\proof By Lemma \ref{out},
$$\int_\Sigma u_k^2dv_g=\|u_k\|_p^2\int_\Sigma \psi_p^2dv_g+o(\|u_k\|_p^2),$$
which leads to
$$\|u_k\|_p^2=\|u_k\|_2^2/(\|\psi_p\|_2^2+o_k(1)).$$
Note that $\|\psi_p\|_2>0$ since $\psi_p\not\equiv 0$. We get the desired result. $\hfill\Box$\\

To proceed, we need gradient estimates for $u_k$, which are analogs of \cite{Druet,Mart,Mar-Stru,L-R-S,Yang-Calc}.
The difference is that $u_k$ changes sign in our case.

\subsection{Gradient estimates}
 Recalling $x_k\ra x_0$ as $k\ra\infty$, we first prove a weak gradient estimate for $u_k$
as below.

\begin{lemma}\label{Prop1}
There exists a constant $C$ depending only on $(\Sigma,g)$ such that for all $x\in \Sigma$ and all $k$,
\be\label{prop-1}\f{\beta_k}{\lambda_k}u_k^2(x)\exp{(\sigma_ku_k^2(x))}({\rm dist}_g(x,x_k))^2\leq C.\ee
\end{lemma}

\noindent{\it Proof.} Suppose not. There would exist $(x_{1,k})\subset \Sigma$ such that
\bea\nonumber\f{\beta_k}{\lambda_k}u_k^2(x_{1,k})\exp{(\sigma_k u_k^2(x_{1,k}))}({\rm dist}_g(x_{1,k},x_k))^2&=&\max_{x\in\Sigma}
\f{\beta_k}{\lambda_k}u_k^2(x)\exp{(\sigma_k u_k^2(x))}({\rm dist}_g(x,x_k))^2\\&\ra&+\infty\quad{\rm as} \quad k\ra\infty.\label{not}\eea
 Let $r_{1,k}>0$ satisfy
\be\label{r-k}r_{1,k}^2\f{\beta_k}{\lambda_k}u_k^2(x_{1,k})\exp{(\sigma_ku_k^2(x_{1,k}))}=1.\ee
It follows easily from (\ref{la-infty}), (\ref{uk-0}), (\ref{not}) and (\ref{r-k}) that $|u_k(x_{1,k})|\ra +\infty$, $r_{1,k}\ra 0$ and $x_{1,k}\ra x_0$. By (\ref{scale}) and
(\ref{r-k}), we have $r_{1,k}\geq r_k$. Also (\ref{not}) leads to
\be\label{compare}\f{{\rm dist}_g(x_{1,k},x_k)}{r_{1,k}}\ra +\infty\ee
as $k\ra\infty$. Take an isothermal coordinate system $(B_{x_{1,k}}(\delta/2), \phi_k; \{y^1,y^2\})$ around $x_{1,k}$, which is constructed
as in Lemma \ref{Lemma-k}, where $x_k$ is replaced by $x_{1,k}$. In particular, $\widetilde{x}_{1,k}=\phi_k(x_{1,k})=0$ and
$\widetilde{x}_k=\phi_k(x_k)\in \phi_k(B_{x_{1,k}}(\delta/2)$.
 Define for $y\in \mathcal{D}_k=\{y\in\mathbb{R}^2: r_{1,k}y\in \phi_k(B_{x_{1,k}}(\delta/2)\}$,
$$v_{1,k}(y)=\f{\widetilde{u}_k(r_{1,k}y)}{u_k(x_{1,k})},$$
where $\widetilde{u}_k=u_k\circ\phi_k^{-1}$. It follows from (\ref{E-L}) and (\ref{lap-nab}) that
\bea\nonumber
-\Delta_{\mathbb{R}^2}v_{1,k}(y)&=&\exp{(f_k(r_{1,k}y))}\le\{\f{v_{1,k}(y)}{u_k^2(x_{1,k})}\exp{(\sigma_k(\widetilde{u}_k^2(r_{1,k}y)
-u_k^2(x_{1,k})))}\ri.\\&&\quad\quad\quad\quad\quad\quad\le.
+r_{1,k}^2\gamma_k v_{1,k}(y)-r_{1,k}^2\f{\mu_k}{u_k(x_{1,k})}\ri\}.\label{vk}
\eea
Let $R$ be any fixed positive number.
Since ${\rm dist}_g(x_{1,k},x_k)=(1+o_k(1))|\widetilde{x}_k|$ and
$|r_{1,k}y-\widetilde{x}_k|=(1+o_k(1)){\rm dist}_g(\phi_k^{-1}(r_{1,k}y),x_k)$
for all $|y|\leq R$,  we have by (\ref{not}) and (\ref{compare}) that
$$
|r_{1,k}y-\widetilde{x}_k|^2\widetilde{u}_k^2(r_{1,k}y)\f{\beta_k}{\lambda_k}\exp{(\sigma_k\widetilde{u}_k^2
(r_{1,k}y))}
\leq (1+o_k(1))|\widetilde{x}_k|^2u_k^2(x_{1,k})\f{\beta_k}{\lambda_k}\exp{(\sigma_k{u}_k^2(x_{1,k}))},
$$
and that for all $|y|\leq R$,
\be\label{ast}v_{1,k}^2(y)\exp{(\sigma_k(\widetilde{u}_k^2(r_{1,k}y)-{u}_k^2(x_{1,k})))}\leq 1+o_k(1).\ee
This leads to
\be\label{stai}\limsup_{k\ra\infty}\max_{|y|\leq R}|v_{1,k}(y)|\leq 1.\ee
Since $R>0$ is arbitrary, having in mind (\ref{stai}),
 we obtain by applying elliptic estimates to (\ref{vk}),
\be\label{v-1k}v_{1,k}\ra 1\quad{\rm in}\quad C^2_{\rm loc}(\mathbb{R}^2).\ee

Define for $y\in\mathcal{D}_k$,
$$\varphi_{1,k}(y)=u_k(x_{1,k})(\widetilde{u}_k(r_{1,k}y)-u_k(x_{1,k})).$$
Then by (\ref{E-L}) and (\ref{lap-nab}),
\bea\nonumber-\Delta_{\mathbb{R}^2}\varphi_k(y)&=&\exp{(f_k(r_{1,k}y))}\le\{v_{1,k}(y)\exp{(\sigma_k(\widetilde{u}_k^2
(r_{1,k}y)-u_k^2(x_{1,k})))}\ri.\\[1.2ex]&&\quad\quad\quad\le.
+r_{1,k}^2\gamma_k u_k^2(x_{1,k})v_{1,k}(y)-r_{1,k}^2\mu_ku_k(x_{1,k})\ri\}.\label{la-phi}\eea
In view of (\ref{ast}) and (\ref{v-1k}), one has
 $\exp{(\sigma_k(2+o_k(1))\varphi_k(y))}\leq 1+o_k(1)$ for all $|y|\leq R$.
 This together with the fact $\sigma_k\ra 4\pi$ implies
\be\label{les-0}\limsup_{k\ra\infty}\varphi_{1,k}(y)\leq 0\quad{\rm uniformly\,\,\,in}\quad y\in\mathbb{B}_0(R).\ee
 It follows from (\ref{v-1k}) that
 \bna
 r_{1,k}^2u_k^2(x_{1,k})&=&\f{r_{1,k}^2}{\pi}\int_{\mathbb{B}_0(1)}u_k^2(x_{1,k})dy\\
 &=&(1+o_k(1))\f{r_{1,k}^2}{\pi}\int_{\mathbb{B}_0(1)}\widetilde{u}_k^2(r_{1,k}y)dy\\
 &=&\f{1+o_k(1)}{\pi}\int_{\mathbb{B}_0{(r_{1,k})}}\widetilde{u}_k^2(y)dy\\
 &\leq&\f{1+o_k(1)}{\pi}\int_\Sigma u_k^2dv_g,
 \ena
 which together with Lemma \ref{concentration} gives
 \be\label{r-u}r_{1,k}^2u_k^2(x_{1,k})=o_k(1).\ee
 Combining (\ref{les-0}) and (\ref{r-u}), we obtain by applying
 elliptic estimates to the equation (\ref{la-phi}),
\be\label{var-con-1}\varphi_{1,k}\ra \varphi\quad{\rm in}\quad C^2_{\rm loc}(\mathbb{R}^2),\ee
where $\varphi$ is defined as in (\ref{bubble}). In view of (\ref{bubble}) and (\ref{var-con-1}),  we have
by noticing that $\beta_k=1+o_k(1)$,
\bna
1+o_k(1)&=&\f{\beta_k}{\lambda_k}\int_\Sigma u_k^2\exp{(\sigma_ku_k^2)}dv_g\\
&\geq&\f{\beta_k}{\lambda_k}\le\{\int_{{B}_{{x}_k}(Rr_k)}{u}_k^2\exp{(\sigma_k{u}_k^2)}dx+
\int_{{B}_{x_{1,k}}(Rr_{1,k})}{u}_k^2\exp{(\sigma_k{u}_k^2)}dx\ri\}\\
&=&(2+o_k(1))\le(\int_{\mathbb{B}_0(R)}\exp{(8\pi\varphi(y))}dy+o_R(1)\ri),
\ena
which is impossible if $R$ is chosen sufficiently large since $\int_{\mathbb{R}^2}\exp{(8\pi\varphi(y))}dy=1$.
This confirms (\ref{prop-1}) and completes the proof of the lemma. $\hfill\Box$\\

We next prove a strong gradient estimate similar to that of DelaTorre-Mancini \cite{D-M}, namely.
\begin{lemma}\label{Prop2}
There exists some constant $C$ such that for all $x\in\Sigma$ and all $k$,
\be\label{gr}{\rm dist}_g(x,x_k)|u_k(x)||\nabla_g u_k(x)|\leq C.\ee
\end{lemma}
\noindent{\it Proof.} Suppose (\ref{gr}) does not hold. Then
$L_k=\sup_{x\in\Sigma}{\rm dist}_g(x,x_k)|u_k(x)||\nabla_g u_k(x)|\ra+\infty$ as $k\ra\infty$. Take $x_{2,k}$ such that
$${\rm dist}_g(x_{2,k},x_k)|u_k(x_{2,k})||\nabla_g u_k(x_{2,k})|=L_k.$$
Noting that $u_k\ra 0$ in $C^2_{\rm loc}(\Sigma\setminus\{x_0\})$, we have $x_{2,k}\ra x_0$ as $k\ra\infty$,
where $x_0$ is given as in (\ref{c-infty}). Take an isothermal coordinate system $(B_{x_{2,k}}(\delta/2),\phi_k;\{y^1,y^2\})$
around $x_{2,k}$, constructed as in Lemma \ref{Lemma-k}. In particular, $\phi_k(x_{2,k})=0$, $\widetilde{x}_k=\phi_k(x_k)\in
\phi_k(B_{x_{2,k}}(\delta/2))\subset\mathbb{R}^2$. In this coordinate system, we have
\be\label{coord}|\widetilde{x}_k||\widetilde{u}_k(0)|
|\nabla_{\mathbb{R}^2}\widetilde{u}_k(0)|=(1+o_k(1))L_k,\ee
where $\widetilde{u}_k=u_k\circ \phi_k^{-1}$.
Denote $s_k=|\widetilde{x}_k|$ and $y_k=\widetilde{x}_k/s_k$.
Assume $y_k\ra \overline{y}$ as $k\ra\infty$, where $|\overline{y}|=1$. Define $v_k(y)=\widetilde{u}_k(s_ky)$
for $y\in \mathcal{U}_k=\{y\in\mathbb{R}^2: s_ky\in \phi_k(B_{x_{2,k}}(\delta/2))\}$. By (\ref{E-L}) and (\ref{lap-nab}),
\bea
\nonumber-\Delta_{\mathbb{R}^2}v_k(y)&=&-s_k^2\Delta_{\mathbb{R}^2}\widetilde{u}_k(s_ky)\\\label{la-v}
&=&s_k^2\exp{(f(s_ky))}\le\{\f{\beta_k}{\lambda_k}\widetilde{u}_k(s_ky)\exp{(\sigma_k
\widetilde{u}_k^2(s_ky))}+\gamma_k\widetilde{u}_k(s_ky)-\mu_k\ri\}.
\eea
By Lemma \ref{Prop1}, there exists a constant $C$ such that for all $y\in\mathcal{U}_k$,
$$|s_ky-\widetilde{x}_k|^2\f{\beta_k}{\lambda_k}\widetilde{u}_k^2(s_ky)
\exp{(\sigma_k\widetilde{u}_k^2(s_ky))}\leq C.$$
This implies that for any fixed $R>0$, there exists some constant $C$ depending on $R$ satisfying
\be\label{ine-1}s_k^2\f{\beta_k}{\lambda_k}\widetilde{u}_k^2(s_ky)
\exp{(\sigma_k\widetilde{u}_k^2(s_ky))}\leq C,\quad\forall\, y\in\mathcal{U}_k\setminus
\mathbb{B}_{\overline{y}}(\f{1}{R}).\ee
Moreover, by a change of variables
\bea\label{ine-2}
\int_{\mathbb{B}_0(R)}|s_k^2\widetilde{u}_k(s_ky)|^pdy\leq (1+o_k(1))s_k^{2p-2}\int_\Sigma|u_k|^pdv_g=o_k(1).
\eea
Also, we have $|\nabla_{\mathbb{R}^2}v_k^2(y)|\leq CL_k$ in $\mathcal{U}_k\setminus
\mathbb{B}_{\overline{y}}(\f{1}{R})$. Hence for all $y\in \mathbb{B}_0(R)\setminus
\mathbb{B}_{\overline{y}}(\f{1}{R})$
\bea
\nonumber v_k^2(y)&\leq&v_k^2(0)+CL_k|y|\\
&\leq&CM_k^2,\label{ine-3}
\eea
where $M_k=|v_k(0)|+|\nabla v_k(0)|$, since by (\ref{coord}), $L_k=(1+o_k(1))|v_k(0)||\nabla_{\mathbb{R}^2}v_k(0)|\leq (1+o_k(1))M_k^2$.
Obviously $M_k\ra\infty$ as $k\ra\infty$.
Define $\overline{v}_k(y)=v_k(y)/M_k$.
In view of (\ref{la-v})-(\ref{ine-3}) and the fact that $\mu_k$ is bounded, we  conclude by using elliptic estimates that $\overline{v}_k\ra v$
in $C^1_{\rm loc}(\mathbb{R}^2\setminus\{\overline{y}\})$, where $v$ is a harmonic function in $\mathbb{R}^2\setminus\{\overline{y}\}$.
Since for any $R>0$,
\bea\nonumber
\int_{\mathbb{B}_0(R)}|\nabla_{\mathbb{R}^2}\overline{v}_k|^2dy&=&\f{1}{M_k^2}\int_{\mathbb{B}_0({Rs_k})}
|\nabla_{\mathbb{R}^2}\widetilde{u}_k(z)|^2dz\\\nonumber
&\leq&\f{1}{M_k^2}(1+o_k(1))\int_\Sigma|\nabla_gu_k|^2dv_g\\\label{xing}
&=&o_k(1),
\eea
we obtain
$$\int_{\mathbb{B}_0(R)\setminus\mathbb{B}_{\overline{y}}(\f{1}{R})}
|\nabla_{\mathbb{R}^2}{v}|^2dy=\lim_{k\ra\infty}\int_{\mathbb{B}_0(R)\setminus\mathbb{B}_{\overline{y}}(\f{1}{R})}
|\nabla_{\mathbb{R}^2}\overline{v}_k|^2dy
=0.$$
Hence $v\equiv C$ for some constant $C$ and $|\nabla_{\mathbb{R}^2}\overline{v}_k(0)|\ra |\nabla_{\mathbb{R}^2}v(0)|=0$. Thus
$|\nabla_{\mathbb{R}^2}{v}_k(0)|=o(M_k)$ and
\be\label{tend-1}\overline{v}_k\ra 1\quad{\rm in}\quad C^1_{\rm loc}(\mathbb{R}^2\setminus\{\overline{y}\}).\ee
Since $L_k=(1+o_k(1))|v_k(0)||\nabla_{\mathbb{R}^2}v_k(0)|\ra+\infty$ as $k\ra\infty$, we have $|\nabla_{\mathbb{R}^2}v_k(0)|>0$ for sufficiently large $k$.
Define  for $y\in \mathcal{U}_k$,
$$v_k^\ast(y)=\f{v_k(y)-v_k(0)}{|\nabla_{\mathbb{R}^2}v_k(0)|}.$$
By (\ref{la-v}), we have
\bna-\Delta_{\mathbb{R}^2}v^\ast_k(y)
=(1+o_k(1))\f{|v_k(0)|}{L_k}s_k^2\exp{(f_k(s_ky))}\le\{\f{\beta_k}{\lambda_k}v_k(y)\exp{(\sigma_k
v_k^2(y))}+\gamma_k v_k(y)-\mu_k\ri\}.\ena
On one hand, by definition of $L_k$ and (\ref{tend-1}), we have
$v^\ast_k\in L^\infty_{\rm loc}(\mathbb{R}^2\setminus\{\overline{y}\})$. On the other hand, by $L_k\ra+\infty$, analogs of (\ref{ine-1}) and (\ref{ine-2}),
we conclude that $\Delta_{\mathbb{R}^2}v^\ast_k(y)\ra 0$ in $L^p_{\rm loc}(\mathbb{R}^2\setminus\{\overline{y}\})$ for any $p>1$.
Then elliptic estimates leads to
$v_k^\ast\ra v^\ast$ in $C^1_{\rm loc}(\mathbb{R}^2\setminus\{\overline{y}\})$, where $v^\ast$ is a harmonic function.
Similarly as (\ref{xing}), we have
$$\int_{\mathbb{B}_0(R)\setminus\mathbb{B}_{\overline{y}}(\f{1}{R})}|\nabla_{\mathbb{R}^2}v^\ast(y)|^2dy\leq \lim_{k\ra\infty}
\int_{\mathbb{B}_0(R)\setminus\mathbb{B}_{\overline{y}}(\f{1}{R})}|\nabla_{\mathbb{R}^2}v_k^\ast(y)|^2dy=0.$$
Noting that $v_k^\ast(0)=0$, we conclude $v^\ast\equiv 0$ in $\mathbb{R}^2\setminus\{\overline{y}\}$. This contradicts
$|\nabla_{\mathbb{R}^2}v^\ast(0)|=1$ and completes the proof of the lemma. $\hfill\Box$

\subsection{Local blow-up analysis}\label{Sec3.4}
  We now consider the local behavior of $u_k$ near $x_0$. In the isothermal coordinate system
  $(B_{x_k}(\delta/2),\phi_k;\{y^1,y^2\})$ defined as in Lemma \ref{Lemma-k}, for convenience, we rewrite (\ref{bubble}) by
  $$\sigma_kc_k(\widetilde{u}_k(2\sigma_k^{-\f{1}{2}}r_ky)-c_k)\ra\varphi_0(y)\quad
{\rm in}\quad C^2_{\rm loc}(\mathbb{R}^2),$$
  where $\varphi_0(y)=-\log(1+|y|^2)$ satisfies
  $$ -\Delta_{\mathbb{R}^2}\varphi_0=4\exp{(-2\varphi_0)}\quad{\rm in}\quad\mathbb{R}^2.$$
  We define
  \be\label{varphi}\varphi_k(y)=\varphi_0\le(\f{\sqrt{\sigma_k}}{2{r}_k}y\ri),\quad y\in\mathbb{R}^2.\ee
  For any fixed $\tau$ with $0<\tau<1$, we let $r_{k,\tau}>0$ be such that
  $\varphi_k(r_{k,\tau})=-\tau\sigma_k c_k^2$, which leads to
  \be\label{rktau} r_{k,\tau}^2=\f{4}{\sigma_k}r_k^2\exp{(\tau \sigma_kc_k^2+o_k(1))}.\ee
  In view of (\ref{conv}), one has $r_{k,\tau}\ra 0$ as $k\ra\infty$.
  Here and in the sequel, we slightly abuse a notation. If $u$ is a function radially symmetric with respect to $0$, we write $u(r)=u(x)$
  with $|x|=r$.

  Let $S_k$ be the radially symmetric solution of
  $$\label{Sk}\le\{\begin{array}{lll}
  -\Delta_{\mathbb{R}^2}S_k(y)=\f{\beta_k}{\lambda_k} S_k(y)\exp{(S_k^2(y))}+\gamma_k S_k(y)-\sqrt{\sigma_k}\,\mu_k\\[1.2ex]
  S_k(0)=\sqrt{\sigma_k}\,{c}_k
  \end{array}\ri.$$
  and $\eta_0$ be the radially symmetric solution of
  $$\le\{\begin{array}{lll}-\Delta_{\mathbb{R}^2}\eta_0=8\eta_0\exp{(2\varphi_0)}+4(\varphi_0^2+\varphi_0)\exp{(2\varphi_0)}\\[1.2ex]
   \eta_0(0)=0.\end{array}\ri.$$
  According to [15], $\eta_0$ satisfies
  \be\label{asy}\eta_0(r)=-\log r^2+O(1)\quad{\rm as}\quad r\ra\infty\ee
   and
  \be\label{whole}-\int_{\mathbb{R}^2}\Delta_{\mathbb{R}^2} \eta_0(y) dy=4\pi.\ee
    Set
  \be\label{eta-k}\eta_k(y)=\eta_0\le(\f{\sqrt{\sigma_k}}{2r_k}y\ri),\quad y\in\mathbb{R}^2.\ee

  For the decomposition of $S_k$, we have the following:
  \begin{lemma}\label{decomp}
  For any $z_k\in \mathbb{B}_{0}(r_{k,\tau})$, there holds
  \be\label{S-decomp}S_k(z_k)=\sqrt{\sigma_k}\,{c}_k+\f{\varphi_k(z_k)}{\sqrt{\sigma_k}\,{c}_k}+\f{\eta_k(z_k)}{\sqrt{\sigma_k}^3{c}_k^3}
  +O\le(\f{1-\varphi_k(z_k)}{{c}_k^5}\ri).\ee
  \end{lemma}
  \noindent{\it Proof.} Let $w_{1,k}$ be given by
  \be\label{S-assume}S_k=\sqrt{\sigma_k}\,c_k+\f{\varphi_k}{\sqrt{\sigma_k}\,{c}_k}+\f{w_{1,k}}{
  \sqrt{\sigma_k}^3{c}_k^3}\ee
  and define $\rho_{1,k}$ by
  \be\label{sca-ineq}\rho_{1,k}=\sup\le\{r: |\eta_k-w_{1,k}|\leq 1-\varphi_k\,\,{\rm in}\,\,[0,r], r\leq r_{k,\tau}\ri\}.\ee
  By (\ref{varphi}), we have
  \be\label{var-k}-\tau\sigma_kc_k^2\leq\varphi_k(r)\leq 0,\quad\forall r\in[0,r_{k,\tau}].\ee
  Also, by (\ref{rktau}), (\ref{asy}) and (\ref{eta-k}), we get
  $$\eta_k(r)=O(c_k^2)\,\,\,{\rm uniformly\,\,in}\,\,\,r\in [0,r_{k,\tau}],$$
  which together with (\ref{sca-ineq}) and (\ref{var-k}) leads to
  \be\label{w-1k}w_{1,k}(r)=O(c_k^2)\,\,\,{\rm uniformly\,\,in}\,\,\,r\in [0,\rho_{1,k}].\ee
   By  (\ref{S-assume}), (\ref{var-k}) and (\ref{w-1k}), we have on $\mathbb{B}_{0}(\rho_{1,k})$ that
  $$\label{geq}S_k\geq (1-\tau)\sqrt{\sigma_k}\,c_k+o_k(1),$$
  that
  \be\label{equ-1}S_k=\sqrt{\sigma_k}\,{c}_k+\f{\varphi_k}{\sqrt{\sigma_k}\,{c}_k}+O\le(\f{1-\varphi_k}{{c}_k^3}\ri),\ee
  and that
  \be\label{geq-2}S_k^2=\sigma_k{c}_k^2+2\varphi_k+\f{\varphi_k^2+2w_{1,k}}{\sigma_k{c}_k^2}+O\le(\f{1+\varphi_k^2}{{c}_k^4}\ri).\ee
  Moreover, we have
  \bna
  -\Delta_{\mathbb{R}^2}w_{1,k}(y)=\sqrt{\sigma_k}^3{c}_k^3\le(\f{\beta_k}{\lambda_k} S_k(y)\exp{(S_k^2(y))}+\gamma_k S_k(y)-
  \sqrt{\sigma_k}\mu_k\ri)+\f{\sigma_k^2c_k^2}{r_k^2}\exp{(2\varphi_k(y))}.
  \ena
  Then using the same argument as the proof of (\cite{Mancini-Thizy}, Step 3.2), and noting that Lemma \ref{Lambda-k}
  is an obvious substitution for (\cite{Mancini-Thizy}, (3.16)), we obtain
  $$\sup_{[0,\rho_{1,k}]}|w_{1,k}-\eta_k|=O\le(\f{1-\varphi_k}{c_k^2}\ri).$$
  This together with (\ref{sca-ineq}) implies $\rho_{1,k}=r_{k,\tau}$ and yields the desired result.
   $\hfill\Box$\\

  Let
  \be\label{rho-k}\rho_k=\sup\le\{r:\widetilde{u}_k(y)\geq \f{1-\tau}{2}c_k,\,{\rm for\,\,all}\,\, |y|\leq r\ri\}.\ee
  We {\it claim} that up to a subsequence
  \be\label{rhok}\rho_k\geq r_{k,\tau}.\ee
  For otherwise, we suppose for all $k$,
  \be\label{rok}\rho_k<r_{k,\tau}.\ee
  By (\ref{scale}) and (\ref{bubble}), there holds
  \be\label{ro-geq}\lim_{k\ra\infty}\f{\rho_k}{r_k}=+\infty.\ee
  Clearly we have for all $y\in \mathbb{B}_0(\rho_k)$
  $$\label{eqution}-\Delta_{\mathbb{R}^2}\widetilde{u}_k(y)=
  \exp{(f_k(y))}\le(\f{\beta_k}{\lambda_k}\widetilde{u}_k(y)\exp{(\sigma_k\widetilde{u}_k^2(y))}+\gamma_k\widetilde{u}_k(y)-\mu_k\ri).$$
  In view of (\ref{bubble}), (\ref{rho-k}), (\ref{rok}), (\ref{ro-geq}), Lemma \ref{Prop2} and
  the fact that $f_k(0)=0$, using an argument of
  (\cite{Druet-Thizy}, Proposition 3.1), we obtain
  \be\label{diff}|\sqrt{\sigma_k} \widetilde{u}_k(y)-S_k(y)|\leq C\f{|y|}{c_k\rho_k} \,\,{\rm for\,\,all}\,\,y\in\mathbb{B}_0(\rho_k).\ee
  By Lemma \ref{decomp}, $S_k(y)\geq (1-\tau+o_k(1))\sqrt{\sigma_k}\,c_k$ for all $y\in\mathbb{B}_0(r_{k,\tau})$. This together with (\ref{diff}) leads to
  $$\widetilde{u}_k(y)\geq \f{2(1-\tau)}{3}c_k\,\,{\rm for\,\,all}\,\,y\in\mathbb{B}_0(\rho_k),$$
  which contradicts (\ref{rho-k}), the definition of $\rho_k$. Hence our claim (\ref{rhok}) holds. Using again the argument
  of (\cite{Druet-Thizy}, Proposition 3.1), we conclude
  \be\label{diff-r}|\sqrt{\sigma_k} \widetilde{u}_k(y)-S_k(y)|\leq C\f{|y|}{c_kr_{k,\tau}} \,\,{\rm for\,\,all}\,\,y\in\mathbb{B}_0(r_{k,\tau}).\ee

  In view of (\ref{varphi}), (\ref{asy}) and (\ref{eta-k}), we have that
  $$\varphi_k=O(1-\eta_k)\quad{\rm in}\quad \mathbb{B}_{0}(r_{k,\tau}).$$
  We have by applying the inequality $|\exp(t)-1-t|\leq t^2\exp{(|t|)}$ for all $t\in\mathbb{R}$ that
  $$\exp\le(\f{\varphi_k^2+2\eta_k}{\sigma_kc_k^2}+O\le(\f{1+\varphi_k^2}{c_k^4}\ri)\ri)=1+\f{\varphi_k^2+2\eta_k}{\sigma_kc_k^2}
  +O\le(\f{(1+\varphi_k^4)\exp(\varphi_k^2/c_k^2)}{c_k^4}\ri),$$
   by employing (\ref{equ-1}) and (\ref{geq-2}) that
  \bna
  \f{\beta_k}{\lambda_k}S_k\exp{(S_k^2)}
  &=&\f{\beta_k}{\lambda_k}\le(\sqrt{\sigma_k}c_k+\f{\varphi_k}{\sqrt{\sigma_k}c_k}+O\le(\f{1-\varphi_k}{c_k^3}\ri)\ri)\\
  &&\times\exp\le(\sigma_k c_k^2+2\varphi_k+\f{\varphi_k^2+2\eta_k}{\sigma_k c_k^2}+O\le(\f{1+\varphi_k^2}{c_k^4}\ri)\ri)\\
  &=&\sqrt{\sigma_k}c_k\f{\beta_k}{\lambda_k}\exp{(\sigma_kc_k^2+2\varphi_k)}\le(1+\f{\varphi_k^2+\varphi_k+2\eta_k}{\sigma_k c_k^2}
  \ri.\\&&\qquad\le.+O\le(\f{(1+\varphi_k^4)\exp(\varphi_k^2/c_k^2)}{c_k^4}\ri)\ri),
  \ena
  and by using (\ref{diff-r}) that
  \bna
  \f{\beta_k}{\lambda_k}\widetilde{u}_k(y)\exp(\sigma_k\widetilde{u}_k^2(y))&=&\f{\beta_k}{\lambda_k}\f{1}{\sqrt{\sigma_k}}
  S_k(y)\exp(S_k^2(y))\le(1+O\le(\f{|y|}{r_{k,\tau}}\ri)\ri)\\
  &=&c_k\f{\beta_k}{\lambda_k}\exp{(\sigma_kc_k^2+2\varphi_k)}\le(1+\f{\varphi_k^2+\varphi_k+2\eta_k}{\sigma_k c_k^2}\ri.\\
  &&\le.+O\le(\f{(1+\varphi_k^4)\exp(\varphi_k^2/c_k^2)}{c_k^4}+\le(1+\f{\varphi_k^2}{c_k^2}\ri)\f{|y|}{r_{k,\tau}}\ri)\ri).
  \ena
  Since $\varphi_k/c_k^2\geq -\tau$ on $\mathbb{B}_0(r_{k,\tau})$, if $\vartheta$ is a real number satisfying $1<\vartheta<2-\tau$, then
  $$(1+\varphi_k^4)\exp{(2\varphi_k+\varphi_k^2/c_k^2)}\leq \exp(\vartheta \varphi_k).$$
  Hence we have by (\ref{scale}) that
  \bna
  \f{\beta_k}{\lambda_k}\widetilde{u}_k\exp(\sigma_k\widetilde{u}_k^2)=\f{\exp{(2\varphi_k)}}{r_k^2c_k}
  \le(1+\f{\varphi_k^2+\varphi_k+2\eta_k}{\sigma_k c_k^2}+
  O\le(\le(1+\f{\varphi_k^2}{c_k^2}\ri)\f{|y|}{r_{k,\tau}}\ri)\ri)+O\le(\f{\exp(\vartheta\varphi_k)}{r_k^2c_k^5}\ri).\label{11}
  \ena
  In the same way, we calculate on $\mathbb{B}_0(r_{k,\tau})$,
  \be
  \f{\beta_k}{\lambda_k}\widetilde{u}_k^2\exp(\sigma_k\widetilde{u}_k^2)
  =\f{\exp{(2\varphi_k)}}{r_k^2}
  \le(1+\f{\varphi_k^2+2\varphi_k+2\eta_k}{\sigma_k c_k^2}+
  O\le(\le(1+\f{\varphi_k^2}{c_k^2}\ri)\f{|y|}{r_{k,\tau}}\ri)\ri)+O\le(\f{\exp(\vartheta\varphi_k)}{r_k^2c_k^4}\ri).\label{u-fu}
  \ee
 Note that
 $$\int_{\mathbb{R}^2}
  \exp{(2\varphi_0(z))}\varphi_0(z)dz=-\pi,$$
  which together with (\ref{varphi}) and (\ref{whole}) leads to
 \bea\nonumber
 &&\int_{\mathbb{B}_0(r_{k,\tau})}\f{\exp{(2\varphi_k)}}{r_k^2}
  \le(1+\f{\varphi_k^2+2\varphi_k+2\eta_k}{\sigma_k c_k^2}\ri)dy\\\nonumber&=&\f{4}{\sigma_k}\int_{\mathbb{B}_0(\f{\sqrt{\sigma_k}}{2r_k}r_{k,\tau})}
  \exp{(2\varphi_0(z))}\le(1+\f{\varphi_0^2(z)+2\varphi_0(z)+2\eta_0(z)}{\sigma_kc_k^2}\ri)dz\\\nonumber
  &=&\f{4}{\sigma_k}\int_{\mathbb{R}^2}\exp{(2\varphi_0(z))}\le(1+\f{\varphi_0^2(z)+2\varphi_0(z)+2\eta_0(z)}{\sigma_kc_k^2}\ri)dz
  +O(\f{r_k}{r_{k,\tau}})\\\nonumber
  &=&\f{4}{\sigma_k}\int_{\mathbb{R}^2}\exp{(2\varphi_0(z))}dz-\f{4}{\sigma_k^2c_k^2}\le(\int_{\mathbb{R}^2}\Delta\varphi_0(z)dz-\int_{\mathbb{R}^2}
  \exp{(2\varphi_0(z))}\varphi_0(z)dz\ri)+O(\f{r_k}{r_{k,\tau}})\\\label{local-1}
  &=&\f{4\pi}{\sigma_k}+O(\f{r_k}{r_{k,\tau}}).
 \eea
  Moreover we have that
  \bea\nonumber\int_{\mathbb{B}_0(r_{k,\tau})}\f{\exp{(2\varphi_k)}}{r_k^2}O\le(\le(1+\f{\varphi_k^2}{c_k^2}\ri)\f{|y|}{r_{k,\tau}}\ri)dy&=&
  \f{r_k}{r_{k,\tau}}\int_{\mathbb{B}_0(\f{\sqrt{\sigma_k}}{2r_k}r_{k,\tau})}\exp{(2\varphi_0(z))}O\le(\le(1+\f{\varphi_0^2(z)}{c_k^2}\ri)
  |z|\ri)dz\\\label{local-2}
  &=&O(\f{r_k}{r_{k,\tau}})
  \eea
  and that
  \bea
  \int_{\mathbb{B}_{0}(r_{k,\tau})}O\le(\f{\exp(\vartheta\varphi_k)}{r_k^2c_k^4}\ri)dy&=&\f{1}{c_k^4}
  \int_{\mathbb{B}_0(\f{\sqrt{\sigma_k}}{2r_k}r_{k,\tau})}O(\exp(\vartheta\varphi_0(z)))dz\nonumber\\\label{local-3}
  &=&O(\f{1}{c_k^4}).
  \eea
  Combining (\ref{u-fu})-(\ref{local-3}), we obtain
  $$\label{energy-1}\int_{\mathbb{B}_{0}(r_{k,\tau})}\f{\beta_k}{\lambda_k}
  \widetilde{u}_k^2\exp(\sigma_k\widetilde{u}_k^2)dy=\f{4\pi}{\sigma_k}+O(\f{1}{c_k^4}),$$
  since $r_k/r_{k,\tau}=O(1/c_k^4)$. It then follows that
  \bea\nonumber
  \int_{\phi_k^{-1}(\mathbb{B}_{0}(r_{k,\tau}))}\f{\beta_k}{\lambda_k}
  {u}_k^2\exp(\sigma_k{u}_k^2)dv_g&=&\int_{\mathbb{B}_{0}(r_{k,\tau})}\f{\beta_k}{\lambda_k}
  \widetilde{u}_k^2(y)\exp(\sigma_k\widetilde{u}_k^2(y))\exp{(f_k(y))}dy\\[1.2ex]\nonumber
  &=&(1+O(r_{k,\tau}))\le(\f{4\pi}{\sigma_k}+O(\f{1}{c_k^4})\ri)\\\label{local-energy}
  &=&\f{4\pi}{\sigma_k}+O(\f{1}{c_k^4}).
  \eea
  Here we used the fact $r_{k,\tau}=O(1/c_k^4)$, which is due to (\ref{conv}) and (\ref{rktau}).

\subsection{Energy estimate away from $x_0$}

To estimate the energy of $u_k$ away from $x_0$, we compute the integral
$$\int_{\Sigma\setminus\phi_k^{-1}(\mathbb{B}_0(r_{k,\tau}))}\f{\beta_k}{\lambda_k}u_k^2\exp(\sigma_ku_k^2)dv_g,$$
where $\phi_k$ is given as in the isothermal coordinate system constructed in Lemma \ref{Lemma-k}.
The following observation is very important for this purpose.
\begin{lemma}\label{Lp}
There exist some $p^\ast>1$ and a constant $C>0$ such that
\be\label{Lp'}\int_{\Sigma\setminus\phi_k^{-1}(\mathbb{B}_0{({r_{k,\tau}})})}\exp(p^\ast\sigma_ku_k^2)dv_g\leq C.\ee
\end{lemma}
\proof Define
$$u_k^{\ast}(x)=\le\{\begin{array}{lll}
u_k(x)&{\rm when}& x\in\Sigma\setminus\phi_k^{-1}(\mathbb{B}_0{({r_{k,\tau}})})\\[1.5ex]
\min\{u_k(x),(1-\tau/2)c_k\}&{\rm when}& x\in\phi_k^{-1}(\mathbb{B}_0{({r_{k,\tau}})}).
\end{array}\ri.$$
In view of (\ref{S-decomp}) and (\ref{diff-r}), there exists a constant $\delta_k>0$ such that if
$x\in\phi_k^{-1}(\mathbb{B}_0{({r_{k,\tau}})})$ satisfies ${\rm dist}_g(x,\p\phi_k^{-1}(\mathbb{B}_0{({r_{k,\tau}})}))\leq\delta_k$,
then there holds  $u_k^\ast(x)=u_k(x)$.
In particular, $0\leq\nu_k=u_k-u_k^\ast\in W_0^{1,2}(\phi_k^{-1}(\mathbb{B}_0(r_{k,\tau})))$ and $u_k^\ast\in W^{1,2}(\Sigma,g)$.
We {\it claim} the following:
\be\label{en-2}\int_\Sigma|\nabla_gu_k^\ast|^2dv_g= 1-\f{\tau}{2}+o_k(1).\ee
By Lemma \ref{concentration}, we have for any $q\geq 1$,
\be\label{mean}\int_\Sigma |u_k^\ast|^qdv_g\leq \int_\Sigma |u_k|^qdv_g=o_k(1).\ee
As a consequence, \be\label{mean-uk}\overline{u_k^\ast}=\f{1}{{\rm Vol}_g(\Sigma)}\int_\Sigma u_k^\ast dv_g=o_k(1).\ee
Testing (\ref{E-L}) by $u_k^\ast-\overline{u_k^\ast}$,
we have by (\ref{bubble}) and (\ref{mean}),
\bea\nonumber
\int_\Sigma|\nabla_gu_k^\ast|^2dv_g&=&\int_\Sigma \nabla_gu_k^\ast\nabla_gu_kdv_g
=\int_\Sigma u_k^\ast\Delta_gu_kdv_g\\\nonumber &=&
\int_{u_k\leq (1-\f{\tau}{2})c_k}\f{\beta_k}{\lambda_k}u_k^2\exp(\sigma_ku_k^2)dv_g+
\int_{u_k> (1-\f{\tau}{2})c_k}\f{\beta_k}{\lambda_k}u_k^\ast u_k\exp(\sigma_ku_k^2)dv_g+o_k(1)\\\nonumber
&\geq& \f{\beta_k}{\lambda_k}\int_{\phi_k^{-1}(\mathbb{B}_0(Rr_k))}(1-\f{\tau}{2})c_ku_k\exp(\sigma_ku_k^2)dv_g+o_k(1)\\\nonumber
&=&(1+o_k(1))(1-\f{\tau}{2})\int_{\mathbb{B}_0(R)}\exp{(8\pi\varphi(y))}dy+o_k(1)\\\label{inequality-1}
&=&1-\f{\tau}{2}+o_k(1)+o_R(1),
\eea
where $o_R(1)\ra 0$ as $R\ra+\infty$.
Testing (\ref{E-L}) by $\nu_k-\overline{\nu_k}$, we have in the same way
\bea\nonumber
\int_\Sigma|\nabla_g\nu_k|^2dv_g&=&\int_\Sigma\nabla_g\nu_k\nabla_g u_kdv_g=\int_\Sigma \nu_k\Delta_gu_kdv_g
\\\nonumber&=&\int_{u_k>(1-\f{\tau}{2})c_k}\le(u_k-(1-\f{\tau}{2})c_k\ri)\f{\beta_k}{\lambda_k}u_k\exp(\sigma_ku_k^2)dv_g+o_k(1)
\\\nonumber
&\geq& (1+o_k(1))\f{\tau}{2}\int_{\mathbb{B}_0(R)}\exp{(8\pi\varphi(y))}dy+o_k(1)\\\label{inequality-2}
&=&\f{\tau}{2}+o_k(1)+o_R(1).
\eea
Moreover,
\be\label{=1}\int_\Sigma|\nabla_gu_k^\ast|^2dv_g+\int_\Sigma|\nabla_g\nu_k|^2dv_g=\int_\Sigma |\nabla_gu_k|^2dv_g=1.\ee
Combining (\ref{inequality-1}), (\ref{inequality-2}) and (\ref{=1}), we conclude our claim (\ref{en-2}).

Let $p_1$ be any fixed number satisfying $1<p_1<1/(1-\tau/2)$. Then it follows from (\ref{en-2})
and Fontana's inequality (\ref{Fontana-ineq}) that
$$\int_\Sigma\exp\le(p_1 (u_k^\ast-\overline{u_k^\ast})^2\ri)dv_g\leq C.$$
 In view of (\ref{mean-uk}), one can choose a number $p^\ast$ such that
$1<p^\ast<p_1$ and
$$\int_\Sigma \exp(p^\ast {u_k^\ast}^2)dv_g\leq C.$$
This particularly implies (\ref{Lp'}) and completes the proof of the lemma.
$\hfill\Box$\\

As a consequence of Lemma \ref{Lp}, we obtain by using the H\"older inequality, $(i)$ of Lemma \ref{Lemma-5}, and Lemma \ref{Lemma-6} that
\bea\nonumber
\int_{\Sigma\setminus\phi_k^{-1}(\mathbb{B}_0(r_{k,\tau}))}\f{\beta_k}{\lambda_k}u_k^2\exp(\sigma_ku_k^2)dv_g&\leq& \f{\beta_k}{\lambda_k}
\|u_k\|_{2q}^2\|\exp(\sigma_ku_k^2)\|_{L^{p^\ast}(\Sigma\setminus\phi_k^{-1}(\mathbb{B}_0(r_{k,\tau})))}\\
&\leq& \f{C}{\lambda_k}\|u_k\|_{2q}^2=o(\|u_k\|_{2q}^4),\label{ou-2}
\eea
where $1/q+1/p^\ast=1$.

\subsection{ Completion of the proof of Proposition \ref{proposition}}

Testing the equation (\ref{E-L}) by $u_k$, in view of (\ref{local-energy}), (\ref{ou-2}) and Lemma \ref{Lemma-6},  we have
\bna\int_\Sigma|\nabla_gu_k|^2dv_g&=&\int_\Sigma u_k\Delta_gu_kdv_g\\&=&\int_\Sigma \gamma_ku_k^2dv_g+\int_{\phi_k^{-1}(\mathbb{B}_0(r_{k,\tau}))}\f{\beta_k}{\lambda_k}u_k^2\exp(\sigma_ku_k^2)dv_g\\
&&\quad+
\int_{\Sigma\setminus\phi_k^{-1}(\mathbb{B}_0(r_{k,\tau}))}\f{\beta_k}{\lambda_k}u_k^2\exp(\sigma_ku_k^2)dv_g\\
&=&\f{4\pi}{\sigma_k}+\gamma_k\|u_k\|_2^2+o(\|u_k\|_{2q}^4).\ena
This together with Lemma \ref{2-p} gives (\ref{Prop-3}), as desired. $\hfill\Box$\\

{\bf Acknowledgements}. This work is partly supported by the National Science Foundation of China (Grant Nos.
  11471014 and  11761131002).

\end{document}